\documentclass[conference]{IEEEtran}
\IEEEoverridecommandlockouts
\usepackage{cite}
\usepackage{amsmath,amssymb,amsfonts}
\usepackage{algorithmic}
\usepackage{graphicx}
\usepackage{textcomp}
\usepackage{xcolor}
\usepackage{algorithm}
\usepackage{multirow}
\usepackage{array}
\usepackage{tikz}
\usepackage{soul}
\def\BibTeX{{\rm B\kern-.05em{\sc i\kern-.025em b}\kern-.08em
		T\kern-.1667em\lower.7ex\hbox{E}\kern-.125emX}}
\begin{document}
	
	\title{Vehicle Route Planning using Dynamically Weighted Dijkstra's Algorithm with Traffic Prediction\\
		\thanks{\textsuperscript{1} Vehicle Engineering Department, Jaguar Land Rover India Limited, Bengaluru 560 048, India} 
		\thanks{\textsuperscript{a} Contributed equally}
		\thanks{\textsuperscript{b} Present address: Department of Engineering, University of Cambridge, Cambridge CB2 1PZ, UK}
		\thanks{\textsuperscript{*} Corresponding author (e-mail: pudhan1@jaguarlandrover.com)}
	}
		\author{\IEEEauthorblockN{Piyush Udhan\textsuperscript{1,a,*}, Akhilesh Ganeshkar\textsuperscript{1,a}, Poobigan Murugesan\textsuperscript{1,a}, Abhishek Raj Permani\textsuperscript{1,a}, \\Sameep Sanjeeva\textsuperscript{1}, Parth Deshpande\textsuperscript{1,b}}
		%\IEEEauthorblockA{\textit{Vehicle Engineering, Jaguar Land %Rover India}\\
			%\IEEEauthorblockA{pudhan1@jaguarlandrover.com, %aganeshk@jaguarlandrover.com, %pmuruges@jaguarlandrover.com, %apermani@jaguarlandrover.com}
		%}
	}
	\maketitle
	
	\begin{abstract}
		Traditional vehicle routing algorithms do not consider the changing nature of traffic. While implementations of Dijkstra’s algorithm with varying weights exist, the weights are often changed after the outcome of algorithm is executed, which may not always result in the optimal route being chosen. Hence, this paper proposes a novel vehicle routing algorithm that improves upon Dijkstra’s algorithm using a traffic prediction model based on the traffic flow in a road network. Here, Dijkstra’s algorithm is adapted to be dynamic and time dependent using traffic flow theory principles during the planning stage itself. The model provides predicted traffic parameters and travel time across each edge of the road network at every time instant, leading to better routing results. The dynamic algorithm proposed here predicts changes in traffic conditions at each time step of planning to give the optimal forward-looking path. The proposed algorithm is verified by comparing it with conventional Dijkstra’s algorithm on a graph with randomly simulated traffic, and is shown to predict the optimal route better with continuously changing traffic.
	\end{abstract}
	\vspace{\baselineskip}
	
	
	\section{Introduction}
	Traffic congestion is a significant challenge in modern society, characterised by slower speeds, longer trip times, and increased vehicle queueing \cite{caves2004encyclopedia}. In a recent report by TomTom, a leading Global Positioning System (GPS) company, the issue of traffic congestion and the resultant loss of time, energy and money of the commuters in major metro cities around the globe is highlighted and the time lost is estimated to be greater than three-fourths of the actual travel time \cite{ccolak2016understanding}. Traffic congestion costed a total of \$160 billion in the United States of America (USA) due to 6.9 billion extra hours travelled and 3.1 billion additional gallons of fuel purchased in 2014. As per the INRIX Roadway Analytics in 2017, the most congested 25 cities of the U.S.A. are estimated to cost the drivers \$480 billion, over the next 10 years, owing to lost time, wasted fuel, and carbon emitted during congestion. In 2018, the total cost of lost productivity due to congestion was found to be \$87 billion in the U.S.A., highlighting the massive effect of traffic congestion on the global economy \cite{afrin2020survey}.
	
	Municipalities, the world over, keep making heavy investments in building new roads, expanding the existing ones and repairing the damaged parts, although controversies on whether broader and a greater number of roads alleviate congestion persist \cite{braess2005paradox}. Traffic congestion also causes increased vehicular emissions of harmful greenhouse gases such as CO\textsubscript{2}, CO, SO\textsubscript{2}, NO\textsubscript{2} which results in a significant rise of air pollutant exposures of commuters and urban populations due to the increased time spent in traffic \cite{xue2013study}. Continued exposure to such harmful emissions can lead to severe respiratory problems and health risks \cite{zhang2013air}\cite{zuurbier2011respiratory}.
	
	One way to reduce traffic congestion is to enable traffic flow optimisation. This would provide insights into more efficient routing in conjunction to smart roads networks. Better routing would eventually lead to better utilisation of all the roads in the network, supporting optimal traffic flow and enabling a reduction in emission. 
	
	Bellman-Ford \cite{bellman1958routing}, A* \cite{hart1968formal}, Floyd Warshall \cite{floyd1962algorithm} and Dijkstra's \cite{dijkstra1959note} are some of the common vehicle routing algorithms. One of the main advantages of Dijkstra’s algorithm is its considerably low complexity, especially when dealing with a sparse road network. There are quite a few efforts studying the improvement of Dijkstra's algorithm and its practical applications. One such effort attempted to solve the Optimal Route Planning in a Parking Lot based on the Dijkstra's algorithm \cite{xiaoxue2017optimal}. The authors combined an impedance function model with the traditional Dijkstra's algorithm to obtain the dynamic time of route. A balance function between distance and time was adopted as the weight matrix. In another paper, modelling of the optimal route-finding problem was carried out using a combination graph with three cost components namely travel distance, toll costs and road surface conditions \cite{soltani2002path}. An optimal route was then calculated based on a combination of these weights fed to the Dijkstra's in the form of an adjacency matrix. In a related study, the authors optimized Dijkstra’s Algorithm from a control network scale and proposed a shortest path algorithm based on the ellipse, which limits the search node collection in a certain area and greatly narrows the search scale \cite{wei2020shortest}.
	
	Although the above studies tried to improve the efficiency of the Dijkstra algorithm, none of them actually solved the vehicle routing problem dynamically. The traffic density, average speed on a patch of road, travel time and traffic flow keep changing dynamically, i.e., at nearly every instant. There is need for a forward-looking traffic model that can predict the various traffic parameters in a road network at future time-steps and feed those in the form of adjacency matrices to the routing algorithm. The routing algorithm itself needs to be modified beyond the conventional Dijkstra’s algorithm, so that it can work with dynamic inputs and provide an optimal route calculated by taking into consideration the future traffic conditions in the road network.
	
	In this paper, traffic flow theory is used for developing a new routing algorithm based on the traffic model of the road network. Consequently, Dijkstra’s algorithm is used here over the A* algorithm, as it traverses through all the nodes of the network. Fundamental diagrams from traffic flow theory provide a good relation between traffic flow rate, traffic density and average speed, that helps model traffic well. This paper also uses the traffic flow conservation law to understand traffic flow at junctions and a Gaussian probabilistic modelling for external traffic inflows into the road network. Here, travel times are chosen as the weights in the adjacency matrix over which optimisation is carried out, instead of a common metric like distance. This is because travel time is an all-encompassing factor, calculation of which takes into consideration all the traffic parameters in the road journey.
	
	The rest of the paper is structured as follows. Section II follows with the methodology used to modify Dijkstra’s algorithm using traffic prediction. Section III discusses the implementation with the steps involved in generating and analysing the traffic. Section IV compares the proposed algorithm with Dijkstra’s algorithm for various routes obtained from real datasets, with concluding remarks in Section V.

\section{Methodology}

\subsection{Dijkstra’s Algorithm}

Dijkstra’s algorithm is an algorithm used for finding the path with the minimum sum of weights between two nodes in a graph. A road network can be represented as a graph. Here, a ‘node’ represents a junction which connects two or more paths, and an ‘edge’ represents a road connecting the ‘nodes’ \cite{dijkstra1959note}. Figures 1 and 2 below depict how a graph containing nodes and edges is generated from an underlying road network obtained from Google Maps.

\begin{figure}[htbp]
	\begin{minipage}[t]{4cm} 
		\centering 
		\includegraphics[scale=0.35]{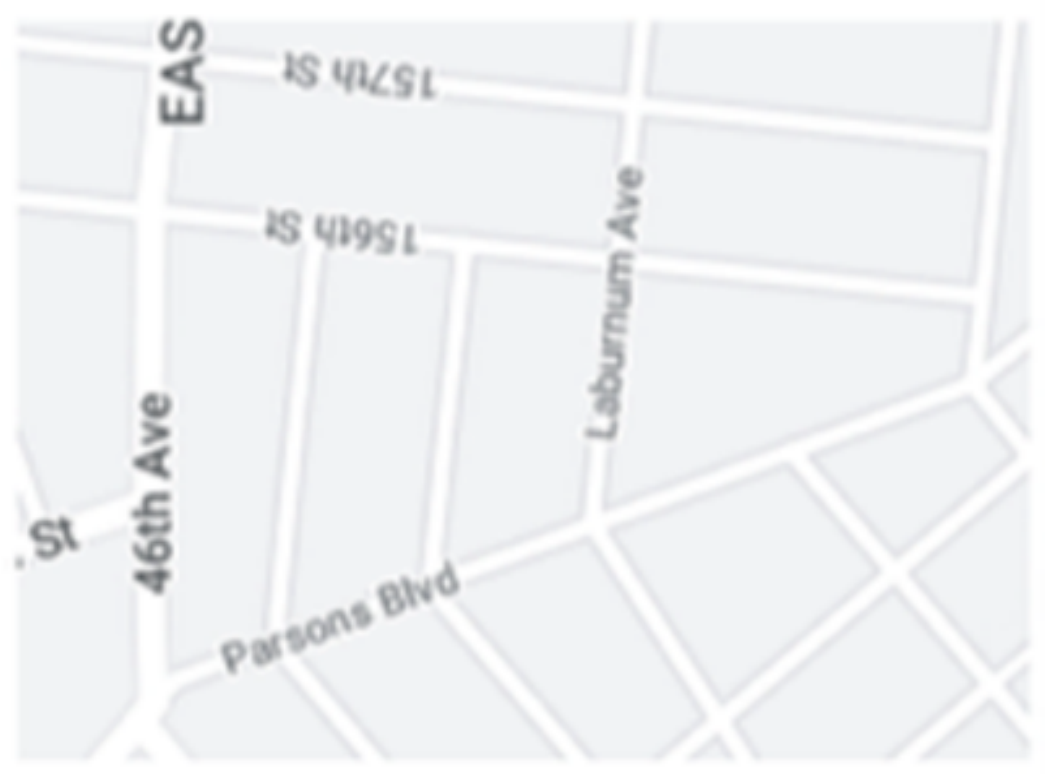} 
		\caption{Actual Road Network} 
	\end{minipage} 
	\hspace{0.3cm} 
	\begin{minipage}[t]{4cm} 
		\centering 
		\includegraphics[scale=0.35]{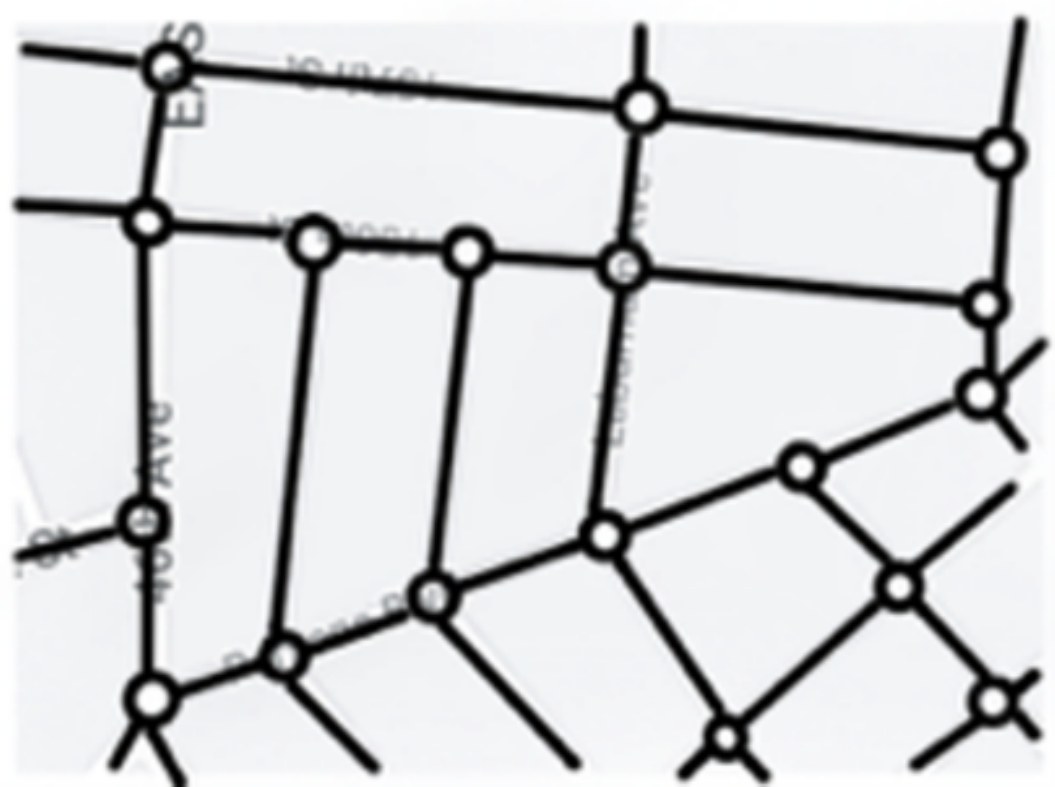} 
		\caption{Derived Graph} 
	\end{minipage} 
\end{figure}

\subsubsection{Modified Dynamic Weight Dijkstra’s Algorithm}

Conventional Dijkstra’s algorithm helps find the shortest path between a source node and a destination node in a static, weighted graph. The algorithm entails visiting each neighbouring node and updating its value, which is a measure of the total weight of the shortest path to it from source, if a smaller distance is found in the iterative process. In this manner, the shortest path from the source node to every other node in a static graph can be found. Considering the graph shown in the Figure 3 as an example, where the source and destination nodes are enumerated as 0 and 9 respectively. The highlighted path, which has been arrived at by visiting neighbouring nodes and iteratively updating their values, is the least weighted route between the two nodes.

\begin{figure}[htbp]
	\centerline{\includegraphics[scale = 0.45]{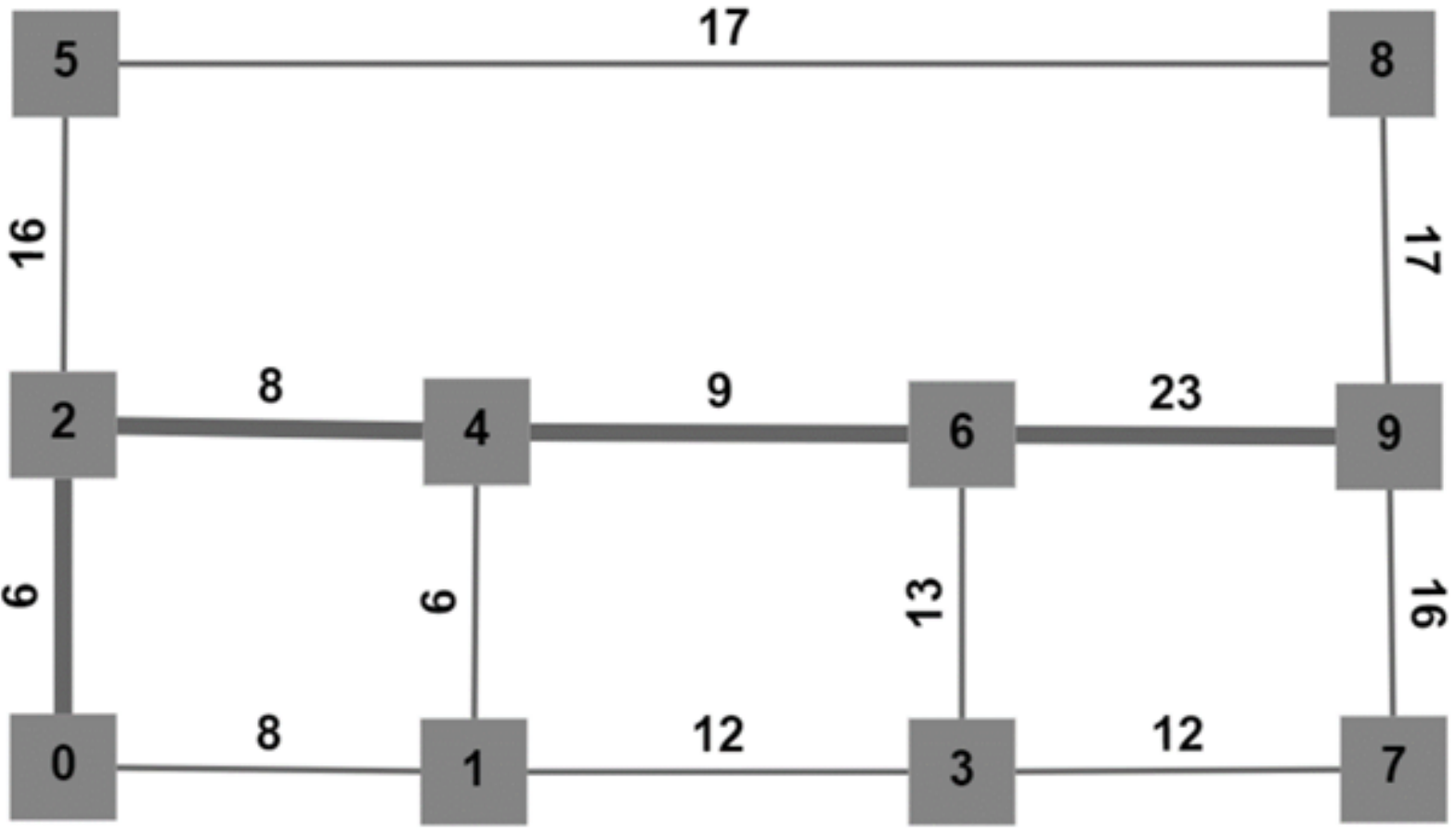}}
	\caption{Shortest Path on Weighted Graph using Dijkstra’s Algorithm}
	\label{fig}
\end{figure}

For the problem statement this paper attempts to solve, the weights of the edges are dynamic and can change with time, in which is representative of a real traffic network. Here, the weight of an edge signifies travel time conforming to traffic flow theory. The inputs to the proposed algorithm are a graph,the source node number and the destination node number. The graph is in the form of a dictionary with the timestamps as keys and an adjacency matrix as the value at each key. The matrix is of N\textsuperscript{th}, where N is the total number of nodes in the graph, and the value of adj[i,j] gives the travel time between nodes i and j at the time instant t as given by the key. This matrix is obtained as an output from the prediction model, which is outlined further in the paper.

The algorithm is an iterative procedure that repeatedly attempts to improve an initial approximation of the values of time[v]. The initial approximation is simply time[0] = 0 (as the source node is indexed as 0) and time[\textit{v}]  = $\infty$ for \textit{v} = 1, 2, …, (\textit{n}-1). In each iteration, a new node is processed, and its time[\textit{v}] value is used to update the time[\textit{v}] value of its immediate successors if the time to reach them is lower through the parent node \textit{v}. Compared to standard Dijkstra’s algorithm, the visited array is removed because this is a dynamic weights problem. So, there might arise a condition where at a later time instant, the weight of an edge falls, requiring the weights of all the connected nodes to be updated, even if they were previously visited. To keep a track of the path taken to reach a particular node minimally, there is a parent array, which is updated with the new parent node if the path to reach the successor is lower through the new node. When the destination node is reached, its neighbours are not pushed back into the queue since we will not leave the destination after reaching it, thus saving computation cycles. The details of this are shown in Algorithm 1.

\begin{algorithm}\caption{Pseudocode for Modified Dynamic Dijkstra}
\begin{algorithmic}[1]
	\renewcommand{\algorithmicrequire}{\textbf{Input:}}
	
	\REQUIRE Graph(dict), source node, destination node
	
	% \STATE first statement
	
	\textit{Initialisation} :
	\FOR {each node \textit{v}}
	\STATE time[\textit{v}] $\leftarrow$ $\infty$
	\STATE parent[\textit{v}] $\leftarrow$ NULL
	\ENDFOR
	\STATE time[source] $\leftarrow$ 0
	\STATE push(0, source) $\to$ \textit{Q}
	\\ LOOP Process
	\WHILE {\textit{Q} \NOT empty}
	\STATE \textit{u} $\leftarrow$ node in \textit{Q} with min time[\textit{u}]
	\STATE remove \textit{u} from \textit{Q}
	\FOR {each neighbour \textit{v} of \textit{u}}
	\STATE \textit{alt} $\leftarrow$ time[\textit{u}] + TRAVELTIME[\textit{u, v}] at time[\textit{u}]
	\IF {\textit{alt} $\textless$ time[\textit{v}]}
	\STATE time[\textit{v}] $\leftarrow$ \textit{alt}
	\STATE parent[\textit{v}] $\leftarrow$ \textit{u}
	\IF {\textit{v} is \textbf{not} destination}
	\STATE push(\textit{alt, v}) $\to$ \textit{Q}
	\ENDIF
	\ENDIF
	\ENDFOR
	\ENDWHILE
	
	\RETURN time[$ $], parent[$ $]
\end{algorithmic}

\end{algorithm}
\label{tab2}

In this case, Dijkstra’s algorithm would use the travel time[\textit{u}, \textit{v}] at time $t$ = 0 which would give an initial path as well as a travel time. However, the traffic parameters would have changed by the time a node \textit{u} is reached. Hence, the travel time at that time instant, which is the output of the prediction model considering updated traffic parameters, is considered in this implementation.

\subsection{Prediction model}
\subsubsection{Traffic Flow Theory}

Traffic flow theory is used to understand various parameters that help define the flow of traffic. The change in traffic parameters of a particular edge over time is attributed to the flow of traffic in and out of the corresponding lane. Each lane is represented as an edge and its weight is defined as the time taken to traverse the lane, which is directly dependent on the traffic in that lane and a few other factors as explained below.

The travel time of the vehicle through a particular lane depends on the length of the lane and the average speed, which in turn depends on the traffic in that lane. Traffic flow theory gives a relationship between the traffic flow rate \textit{q}, density \textit{k} and speed \textit{u}  given by

\begin{equation}
	q(i, j)=k(i, j) \cdot u(i, j)
\end{equation}

where, \textit{q}($i$,$j$) is the traffic flow rate between nodes \textit{i} and \textit{j}, defined as the average number of vehicles passing a specific point per unit time. Similarly, \textit{k}($i$,$j$) is the traffic density, defined as the average number of vehicles occupying a section of the road per unit distance, and \textit{u}($i$,$j$) is the space mean speed defined as the spatial average speed of all vehicles.

The conservation law states that the measurable parameters of an isolated physical system remain constant. In this system, the conserved parameter is the number of vehicles on a given stretch of road. Since the number of vehicles on a given stretch of road between [\textit{x, x'}] is constant, for there to be a change in the traffic density $\rho$, there must be an inflow  or an outflow, such that the traffic flow is conserved, as illustrated in Figure 4.

\begin{figure}[htbp]
	\centerline{\includegraphics[scale = 0.45]{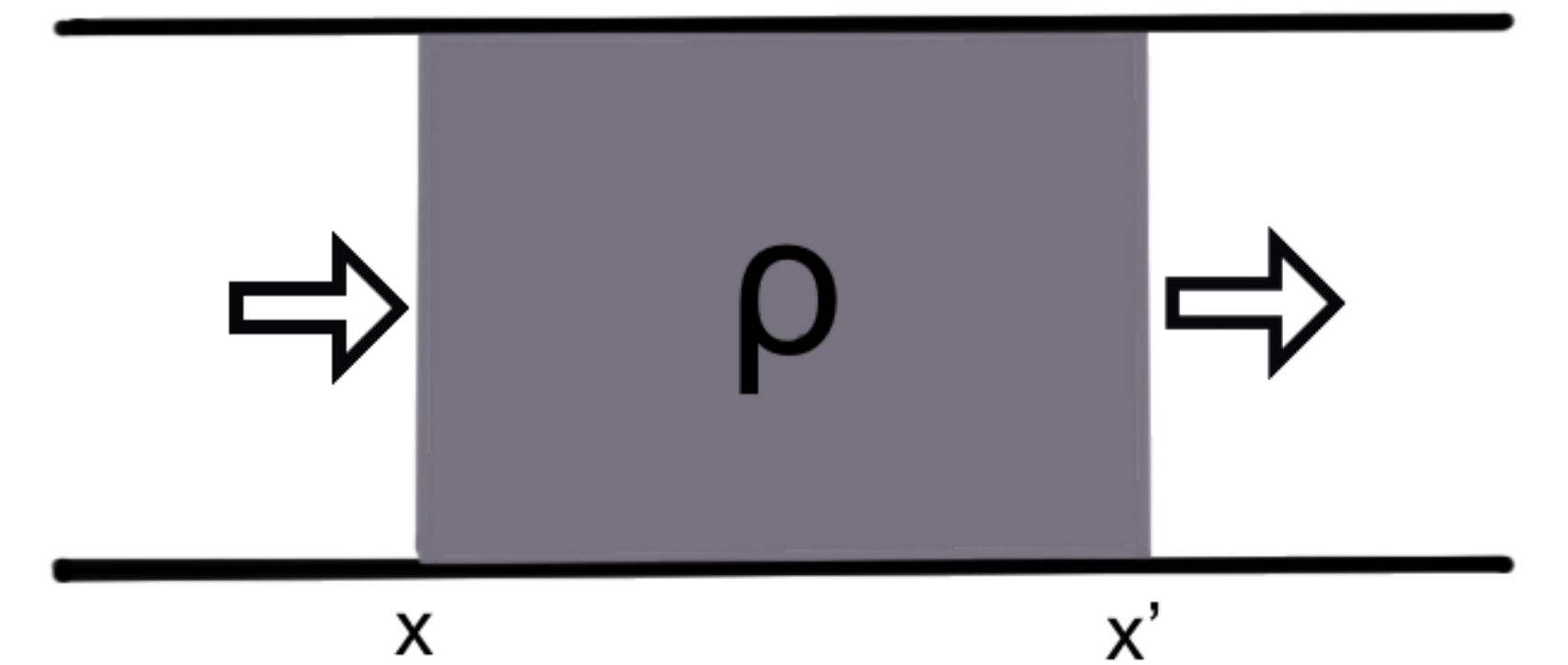}}
	\caption{Conservation law illustration}
	\label{fig}
\end{figure}

In 1955, Lighthill, Whitman and Richards proposed a macroscopic traffic flow model known as the LWR model \cite{lighthill1955kinematic}\cite{ali2016inflow}. According to this model, if a vehicle can be assumed to be a molecule, then the traffic in the lane can be presented as an incompressible fluid which cannot be compressed beyond a certain density. Based on the conservation laws and the relation between the traffic variables, traffic flow is represented by the model as follows

\begin{equation}
	\frac{\partial k}{\partial t}+Q^{\prime}(k) \cdot \frac{\partial k}{\partial x}=f(),
\end{equation}

where, $k$ is the density, $Q$ is the flux or the flow rate, and $f()$ represents the source term which represents the inflow and outflow.

This methodology involves the representation of the system in a state space form. The associated model then uses the conservation law and dynamic speed formula with a linear traffic stream model to study the traffic flow. The density and aggregate space mean speed are considered to be the state variables of interest in this model. The relative difference in the entry and exit flows is taken as the input variable and the space mean speed is the output variable. The classified volume per minute and spot speeds are extracted from the data collected, and the traffic variables, flow and space mean speed, are calculated. The actual densities are calculated using an input-output analysis.

Hence, the traffic flow at the next timestep is determined by considering the present inflows and outflows.

\section{IMPLEMENTATION}

\subsection{Definitions}
\begin{table}[]
	\caption{Definitions}
	\begin{center}
	\begin{tabular}{|l|l|}
		\hline
		\textbf{Parameter}    & \textbf{Definition}                                                                                                                                 \\ \hline
		Edge ID               & \begin{tabular}[c]{@{}l@{}}Unique number assigned to a given edge in \\ the graph for identification purpose.\end{tabular}                           \\ \hline
		Node ID               & \begin{tabular}[c]{@{}l@{}}Unique number assigned to a given node in \\ the graph for identification purpose.\end{tabular}                           \\ \hline
		Edge thickness        & Width of the road represented by given edge.                                                                                                         \\ \hline
		Vehicle thickness     & Width of given vehicle.                                                                                                                              \\ \hline
		Traffic density       & \begin{tabular}[c]{@{}l@{}}Number of vehicles present on a unit length \\ of road.\end{tabular}                                                      \\ \hline
		Edge traffic          & Traffic density on the given edge.                                                                                                                   \\ \hline
		Edge forward traffic  & \begin{tabular}[c]{@{}l@{}}Traffic density calculated when traversing\\ from start node to end node of given edge.\end{tabular}                      \\ \hline
		Edge backward traffic & \begin{tabular}[c]{@{}l@{}}Traffic density calculated when traversing \\ from end node to start node of given edge.\end{tabular}                     \\ \hline
		Edge free flow speed  & \begin{tabular}[c]{@{}l@{}}Max allowed speed of vehicle when given \\ edge is in free-flow state i.e., has very low \\ traffic density.\end{tabular} \\ \hline
		Edge jam speed        & \begin{tabular}[c]{@{}l@{}}Max speed of vehicle when given edge is \\ in jam state i.e., has very high traffic density.\end{tabular}                 \\ \hline
		$\alpha$, $\beta$                & \begin{tabular}[c]{@{}l@{}}Hyper-parameters used for tuning prediction \\ model.\end{tabular}                                                        \\ \hline
	\end{tabular}
\end{center}
\end{table}

\subsection{Prediction Pipeline}
In this section, the architecture of the proposed prediction model is briefly introduced. The prediction model consists of two parts, data simulation and prediction, as shown in Figure 5. The road network is forwarded to the data simulation component in a suitable format, which is then used to generate data at time $t$ = 0 and pass it to the prediction component for predicting traffic flow at future time instances. The components are designed to be independent of each other.

\begin{figure}[htbp]
	\centerline{\includegraphics[scale = 0.8]{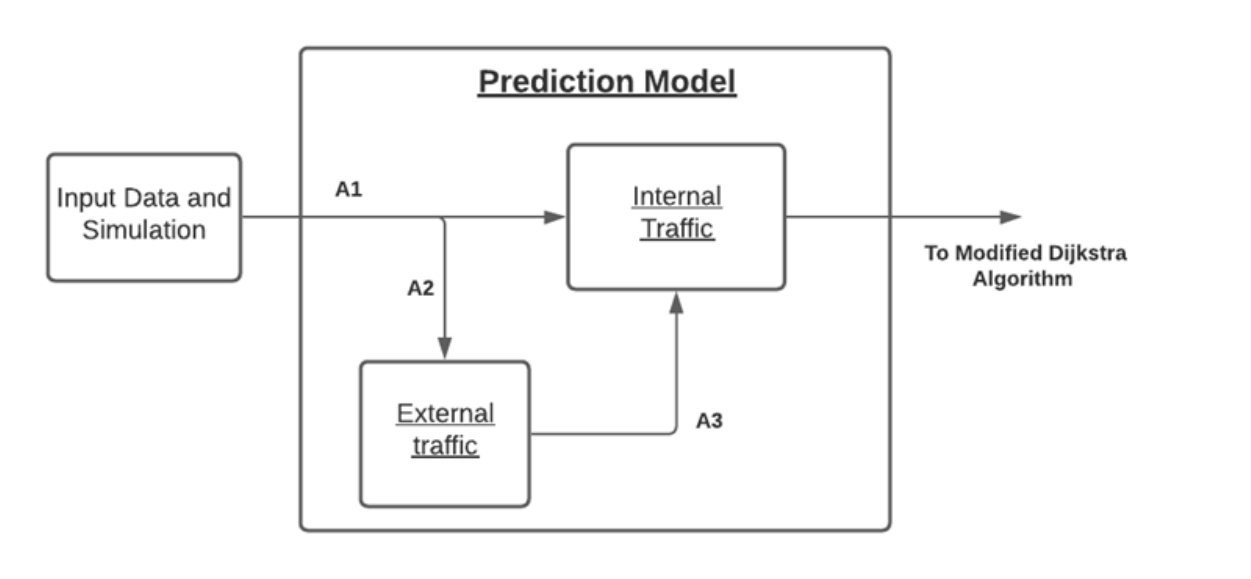}}
	\caption{Prediction Model Pipeline}
	\label{fig}
\end{figure}

\subsection{Input Data and Simulation}
A road network represented in the form of a graph is fed as an input to the model. All the information about the graph is organised into two files: one providing information about the nodes and the other about the edges and edge weights. Table 2 provides details about various classes used in the model and their significance in the model. The classes described as basic requirements are required for expressing the road network in the form of a graph. The table also shows the parameters that are used as outputs.

%\begin{figure}[htbp]
%	\centerline{\includegraphics[scale = 0.65]{Table1.jpg}}
%	\caption{Class Information}
%	\label{Table 1}
%\end{figure}

\begin{table}[]
	\caption{Class Information}
	\begin{tabular}{|l|l|l|}
		\hline
		\multicolumn{1}{|c|}{\textbf{Class}}                    & \multicolumn{1}{c|}{\textbf{Parameters}}                                                                                                 & \multicolumn{1}{c|}{\textbf{Significance}}                                                                                                     \\ \hline
		\multirow{3}{*}{\textbf{Node}}                          & \begin{tabular}[c]{@{}l@{}}Basic: ID, Latitude, \\ Longitude, Category\end{tabular}                                                      & \multirow{2}{*}{Basic requirements}                                                                                                            \\ \cline{2-2}
		& \begin{tabular}[c]{@{}l@{}}Routing: Node connections, \\ Edge connections\end{tabular}                                                   &                                                                                                                                                \\ \cline{2-3} 
		& \begin{tabular}[c]{@{}l@{}}External Traffic: Gaussian mean, \\ Gaussian sigma, Number of \\ external vehicles\end{tabular}               & \begin{tabular}[c]{@{}l@{}}Calculate external \\ traffic\end{tabular}                                                                          \\ \hline
		\multirow{4}{*}{\textbf{Edge}}                          & \begin{tabular}[c]{@{}l@{}}Basic: ID, Start node, End node, \\ Distance, Category\end{tabular}                                           & \multirow{2}{*}{\begin{tabular}[c]{@{}l@{}}Basic requirements\\ for graph\end{tabular}}                                                        \\ \cline{2-2}
		& \begin{tabular}[c]{@{}l@{}}Category: Thickness, \\ Free flow speed, Jam speed\end{tabular}                                               &                                                                                                                                                \\ \cline{2-3} 
		& \begin{tabular}[c]{@{}l@{}}Traffic: Forward traffic, \\ Backward traffic, \\ Vehicle forward list, \\ Vehicle backward list\end{tabular} & \begin{tabular}[c]{@{}l@{}}Required to calculate \\ current traffic flow in \\ a particular edge and \\ to predict future traffic\end{tabular} \\ \cline{2-3} 
		& \begin{tabular}[c]{@{}l@{}}Output Parameters: Forward \\ travel time, Backward travel \\ time\end{tabular}                               & \begin{tabular}[c]{@{}l@{}}Predict internal traffic \\ flow\end{tabular}                                                                       \\ \hline
		\multicolumn{1}{|c|}{\multirow{3}{*}{\textbf{Vehicle}}} & \begin{tabular}[c]{@{}l@{}}Basic: ID, Edge ID, Category, \\ State\end{tabular}                                                           & \multirow{2}{*}{\begin{tabular}[c]{@{}l@{}}Class required to \\ simulate traffic in the \\ given road network\end{tabular}}                    \\ \cline{2-2}
		\multicolumn{1}{|c|}{}                                  & \begin{tabular}[c]{@{}l@{}}Category: Thickness, \\ Maximum speed\end{tabular}                                                            &                                                                                                                                                \\ \cline{2-3} 
		\multicolumn{1}{|c|}{}                                  & \begin{tabular}[c]{@{}l@{}}Prediction: Current speed, \\ Time to complete \\ (current edge/road)\end{tabular}                            & \begin{tabular}[c]{@{}l@{}}Calculate current \\ traffic on a \\ given edge\end{tabular}                                                        \\ \hline
	\end{tabular}
\end{table}

%INSERT TABLE 1 HERE

At time $t$ = 0, the traffic data simulated consists of vehicle positions that are generated randomly over the graph by specifying the number of vehicles as an argument. Vehicles that enter the network from outside it are characterized as external traffic, and the vehicle data generated at every node has an open edge for external traffic inflow. Thus, at time $t$ = 0, the graph and traffic data are provided to the prediction component.

\subsection{Prediction Component}
The task of the prediction component is broken down as predicting the future states of external and internal traffic.\\
\subsubsection{External traffic} 
A Gaussian probability distribution function (GPDF) is used here to simulate external traffic. GPDFs are used to represent the probability distribution function of a normally distributed random variable, with a mean $\mu$ and variance ${\sigma^2}$. 

The Gaussian equation is given by\\
\begin{equation}
	p(x) = \frac{1}{\sigma\sqrt{2\pi}}
	\exp\left( -\frac{1}{2}\left(\frac{x-\mu}{\sigma}\right)^{\!2}\,\right) ,
	\forall x \in X,
\end{equation}
where, $\mu$ is the mean, ${\sigma^2}$ is the variance, and X is the domain of the random variable.

High-order Gaussian process dynamical models for traffic flow prediction \cite{zhao2016high} and a probabilistic model fusing multisource data on road traffic \cite{lin2017road} are two models which use Gaussian models extensively for both prediction and data aspects of their model. These have provided valuable insights and superior results regarding the use of Gaussian models in the vehicle routing problem. Hence, this study uses the Gaussian model.

It is assumed that external traffic enters the road network from a node at time $t$ = 0. The value corresponding to the number of vehicles about to enter the graph from a particular node is used as the mean value for the GPDF, and the variance is initialised in the range of 0 to 3. This is done for every node at which external traffic can enter the road network. After every time instance, this step is repeated, and the number of vehicles generated in the last iteration is assigned as the new mean for the GPDF for that node. Consequently, external traffic is generated at nodes in the graph and passed on to the internal traffic component.

\begin{table}[]
	\caption{Prediction Model}
	\begin{tabular}{|l|l|l|l|}
		\hline
		\textbf{Step} & \textbf{Input}                                                        & \textbf{Principles}                                                     & \textbf{Output}                                                      \\ \hline
		1             & \begin{tabular}[c]{@{}l@{}}Vehicle positions\\ Edge data\end{tabular} & Flow conservation                                                       & Traffic Density                                                      \\ \hline
		2             & \begin{tabular}[c]{@{}l@{}}Traffic density\\ Edge data\end{tabular}   & Fundamental diagrams                                                    & \begin{tabular}[c]{@{}l@{}}Average speed\\ Traffic flow\end{tabular} \\ \hline
		3             & \begin{tabular}[c]{@{}l@{}}Average speed\\ Edge weight\end{tabular}   & \begin{tabular}[c]{@{}l@{}}Distance-speed-time \\ equation\end{tabular} & Travel time                                                          \\ \hline
	\end{tabular}
\end{table}

\subsubsection{Internal Traffic} 
This section covers how internal traffic is predicted using the data provided at time $t$ = 0. The graph (road network) is received via data buses A1 and A3 from the Data and External traffic components, respectively, as shown in Figure 5.
The steps for predicting the internal traffic specified in Table 3 are as follows.

Step 1: Taking vehicle positions and edge data as arguments, the model applies traffic flow conservation. Using the vehicle positions, it calculates the number of vehicles present along every edge. Furthermore, using vehicle characteristics (like thickness) and edge characteristics, the traffic density is calculated at every edge. At this stage, the value assigned to the edge, i.e., the traffic along that edge, is the ratio of the number of vehicles to the edge thickness.

\begin{algorithm}\caption{Determine Speed of Vehicles}
	\begin{algorithmic}[1]
		\renewcommand{\algorithmicrequire}{\textbf{Input:}}
		
		\REQUIRE \text{Traffic\_density}
		
		% \STATE first statement
		
		\textit{Initialisation} :
		\STATE $v_{ret}$ $\leftarrow$ free-flow speed of the edge, where the vehicle is present
		
		\STATE curr\_tr\_dens $\leftarrow$ traffic\_density $\times$ vehicle\_thickness 
		\vspace{\baselineskip}
		
		\STATE curr\_tr\_dens $\leftarrow$ $\dfrac{\textrm{ curr\_tr\_dens}}{\textrm{edge\_length}}$ $\times$ edge\_thickness
		\vspace{\baselineskip}
		
		\IF {curr\_tr\_dens  $\textless$ $\beta$ \& curr\_tr\_dens $\leq$ $\alpha$}
		\STATE $v_{ret}$ $\leftarrow$ (edge.free\_flow\_speed - edge.jam\_speed)
		$\times$ (1 – curr\_tr\_dens) +  edge.jam\_speed
		\ELSIF {curr\_tr\_dens $>$ $\beta$}
		\STATE $v_{ret}$ $\leftarrow$ jam\_speed of edge is assigned
		\ENDIF
		\RETURN $v_{ret}$
		\COMMENT {speed value for vehicle returned}

		\RETURN time$[$ $]$, parent$[$ $]$
	\end{algorithmic}
	
\end{algorithm}
\label{tab2}

Step 2: Using the traffic densities, the speed of every vehicle is calculated using the function shown in Algorithm 2. The average speed is then calculated by taking the mean of the speed of all vehicles along that particular edge.

Step 3: Using the average speed of an edge and the edge length, the travel time for a particular edge is calculated and assigned as the edge-weight for that particular edge.

\subsubsection{Other background processes}
When a vehicle reaches the end of an edge, i.e., it reaches a node $X$, it needs to decide which edge it will travel along next. Hence, the model then analyses the traffic in all the edges directly connected to node $X$, by calling the function described in Algorithm 4 which uses the transformation function defined in Algorithm 3. In Algorithm 3, $\alpha$ and $\beta$ are hyper-parameters used for tuning the prediction model. In doing so, a probability distribution is generated that is used to select the next edge to travel to. A vehicle can also exit the network when it reaches an open node. In such a case, it chooses an edge with edge id = (-1).
Output: For every time instant, the predicted travel time for every edge is assigned as the weight of that edge, and the transformed graph is then passed on to the routing algorithm for route planning.

\begin{algorithm}\caption{Transformation Function Logic}
\begin{algorithmic}[1]
	\renewcommand{\algorithmicrequire}{\textbf{Input:}}
	\REQUIRE Edge traffic
	
	% \STATE first statement
	\COMMENT{two thresholds predefined, $\alpha$ and $\beta$ where $\beta$ $>$ $\alpha$, $\alpha$ and $\beta$ are thresholds for edge states\\
		- traffic less than $\alpha$ $\Rightarrow$ edge in free flow state\\
		- traffic greater than $\alpha$ and less than $\beta$ $\Rightarrow$ edge is in normal flow state\\
		- traffic greater than $\beta$ $\Rightarrow$ edge is in jam state\\
		
		based on state, a quadratic functional diagram transformation is applied on edge traffic and this transformed value is stored in \textit{result}}

	\RETURN \textit{result}

\end{algorithmic}
\end{algorithm}	

\begin{algorithm}\caption{Algorithm for Selecting Next Edge}
\begin{algorithmic}[1]
	\renewcommand{\algorithmicrequire}{\textbf{Input:}}
	\REQUIRE node $X$
	
	% \STATE first statement
	
	\STATE $K$ $\leftarrow$ dictionary with edge id as key and traffic situation as value
	\STATE sum $\leftarrow$ 0
	\FOR {every edge $E$ directly connected to node $X$}
	\STATE tr $\leftarrow$ get value of traffic in that edge
	\IF {node $X$ $\equiv$ start node of edge $E$}
	\STATE tr $\leftarrow$ transformation\_function($E$.forward\_traffic)
	\ELSE
	\STATE tr $\leftarrow$ transformation\_function($E$.backward\_traffic)
	\ENDIF
	\STATE $K$[$E$.id] = tr
	\STATE sum += $K$[$E$.id]
	\ENDFOR
	\COMMENT {Normalising values in $K$}
	\FOR {every edge $E$ id in $K$}
	\STATE $K$[$E$.id]/ = sum
	\ENDFOR
	\COMMENT {Using these values assigned for every directly connected edge in dictionary $K$ as probabilities}
	\STATE choice = random.choices(node\_connections, weights=$k$)
	\COMMENT{this returns chosen edge id}
	\RETURN choice
	
\end{algorithmic}
\end{algorithm}

	\section{EXPERIMENTATION, RESULTS AND OBSERVATIONS}
	In the dataset \cite{sstd2005} used for testing, the road networks of five major cities of the globe are expressed in the form of a graph. The junctions are exhibited as nodes, and the roads connecting them as edges. The rationale behind using this dataset for road networks is to calibrate and evaluate the model on data which mimics real life road networks and traffic conditions as much as possible. As this dataset is designed to replicate actual road networks, it serves as an ideal data feed for the model.
	
	\begin{figure}[htbp]
		\centerline{\includegraphics[scale=0.45]{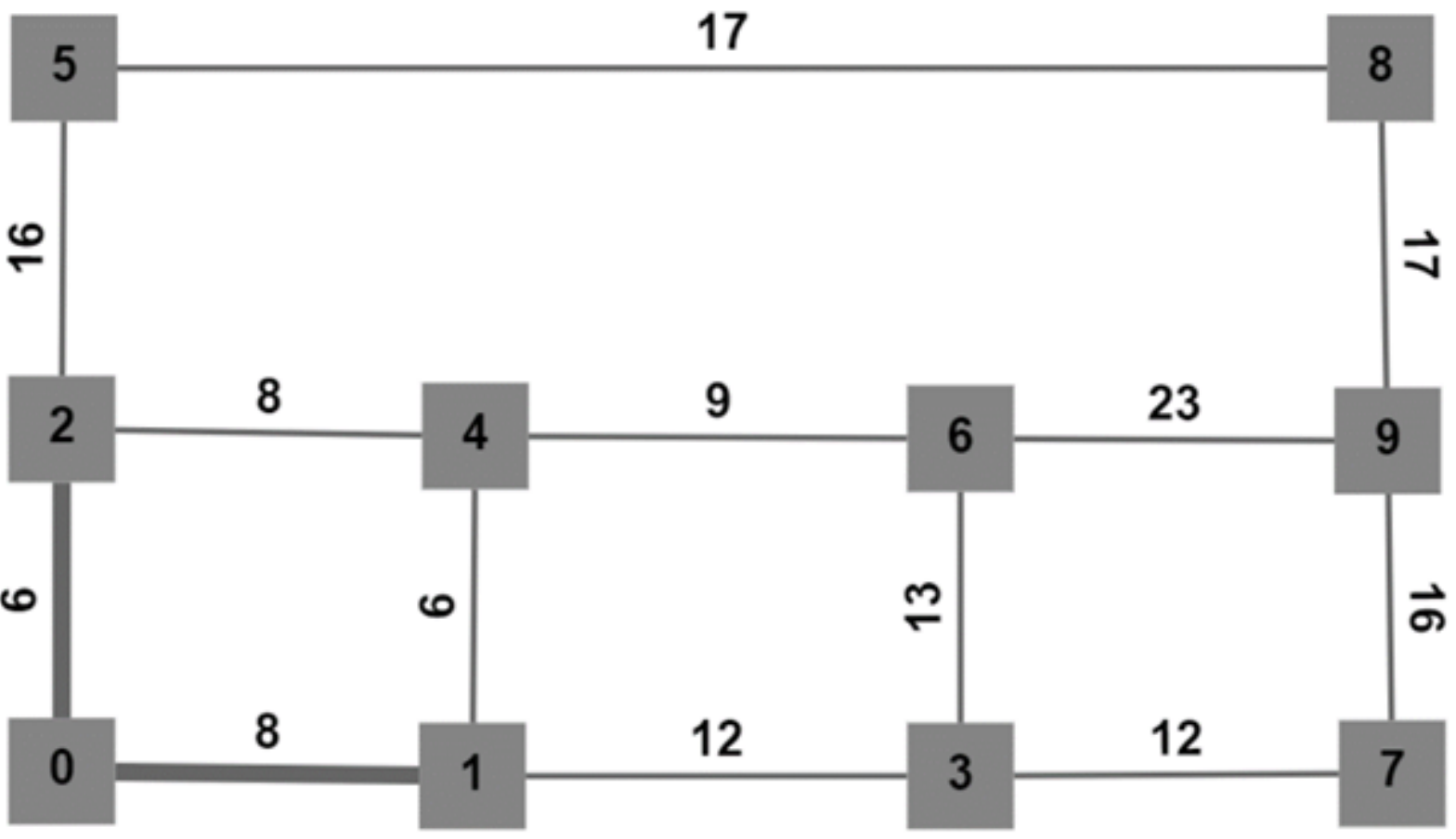}}
		\caption{Road Network at t = 0}
	\end{figure}
	
	As covered in the input data and simulation component of the prediction pipeline in Section 3, random traffic data is generated over the road network. The model is tested with different $\langle \text{source\_node}, \text{destination\_node} \rangle$ pairs to get the total travel time as an output for both, the model proposed in this paper and Dijkstra’s algorithm.
	
	In the case of conventional Dijkstra’s algorithm, the shortest path is computed by applying the unmodified routing algorithm on the graph shown in Figure 7 at t = 0, without considering the dynamic nature of traffic and resulting dynamic edge weights.
	
	Figures 8 and 9 depict the implementation of the proposed routing algorithm. This algorithm considers the dynamic nature of traffic over the road network in subsequent timestamps. It computes the total time required for every route to reach destination node by considering the traffic at future time steps from the prediction model. All the possible edge combinations at different time stamps are computed to obtain the shortest route.

	As seen in the graph, there are two possible paths from the source node 0, E\textsubscript{0-1} and  E\textsubscript{0-2}. The proposed algorithm evaluates each of these two edges at t = 0, and then further evaluates the following edges at corresponding future time stamps based on the dynamic traffic and edge weight inputs from prediction model. The detailed computation data is presented in Table 4.

\begin{figure}[htbp]
	\begin{minipage}[t]{4cm} 
		\centering 
		\includegraphics[scale=0.25]{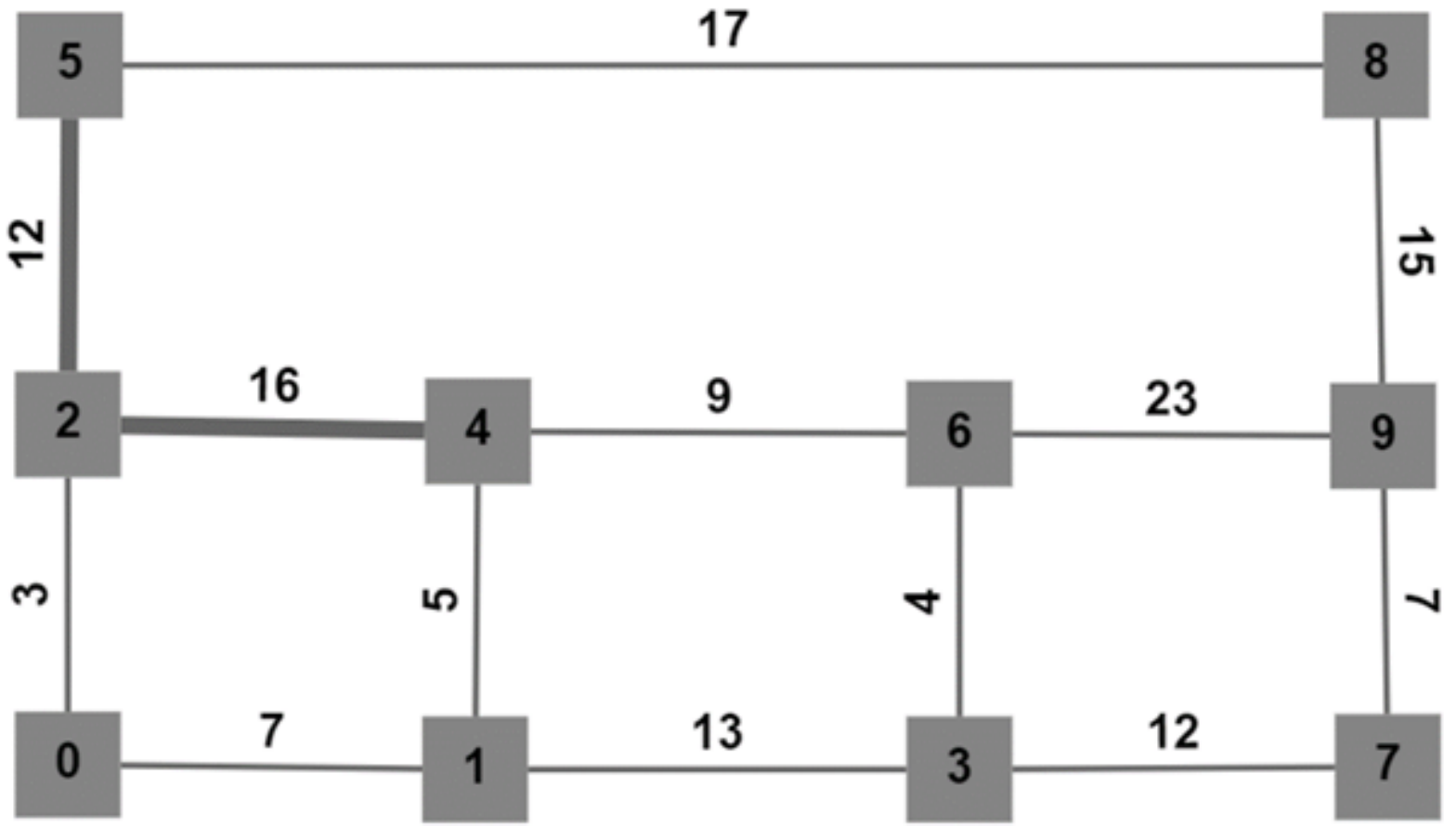} 
		\caption{Road Network at t = 6} 
	\end{minipage} 
	\hspace{0.4cm} 
	\begin{minipage}[t]{4cm} 
		\centering 
		\includegraphics[scale=0.25]{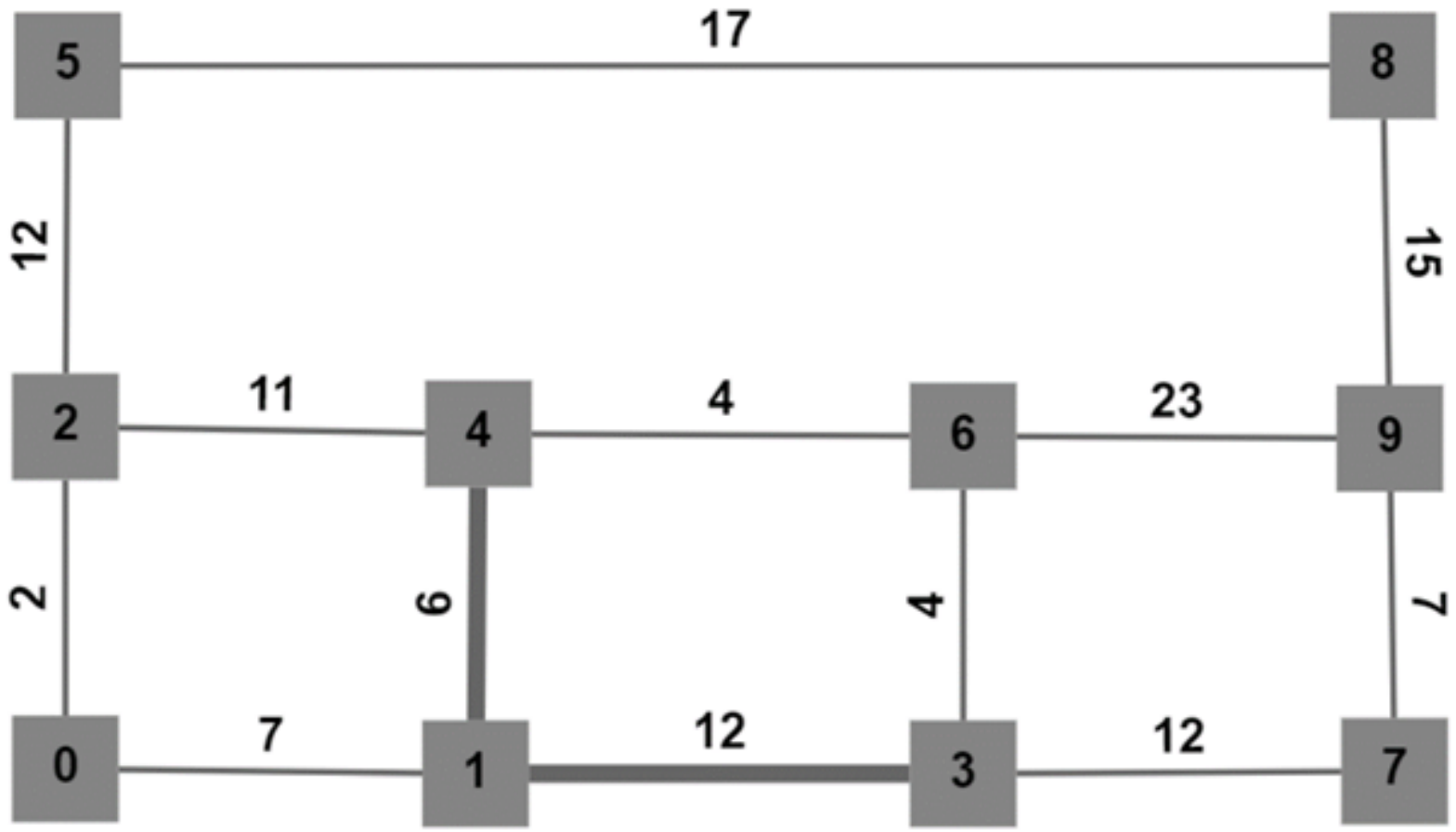} 
		\caption{Road Network at t = 8} 
	\end{minipage} 
\end{figure}
		
	The graph changes at every time instant, as the weights of the edges of the graph change with time. The proposed routing algorithm evaluates all the dynamic weights and time stamp permutations to determine the shortest path from source to destination, as shown in Figures 10 and 11.
	
	\begin{figure}[htbp]
		\begin{minipage}[t]{4cm} 
			\centering 
			\includegraphics[scale=0.25]{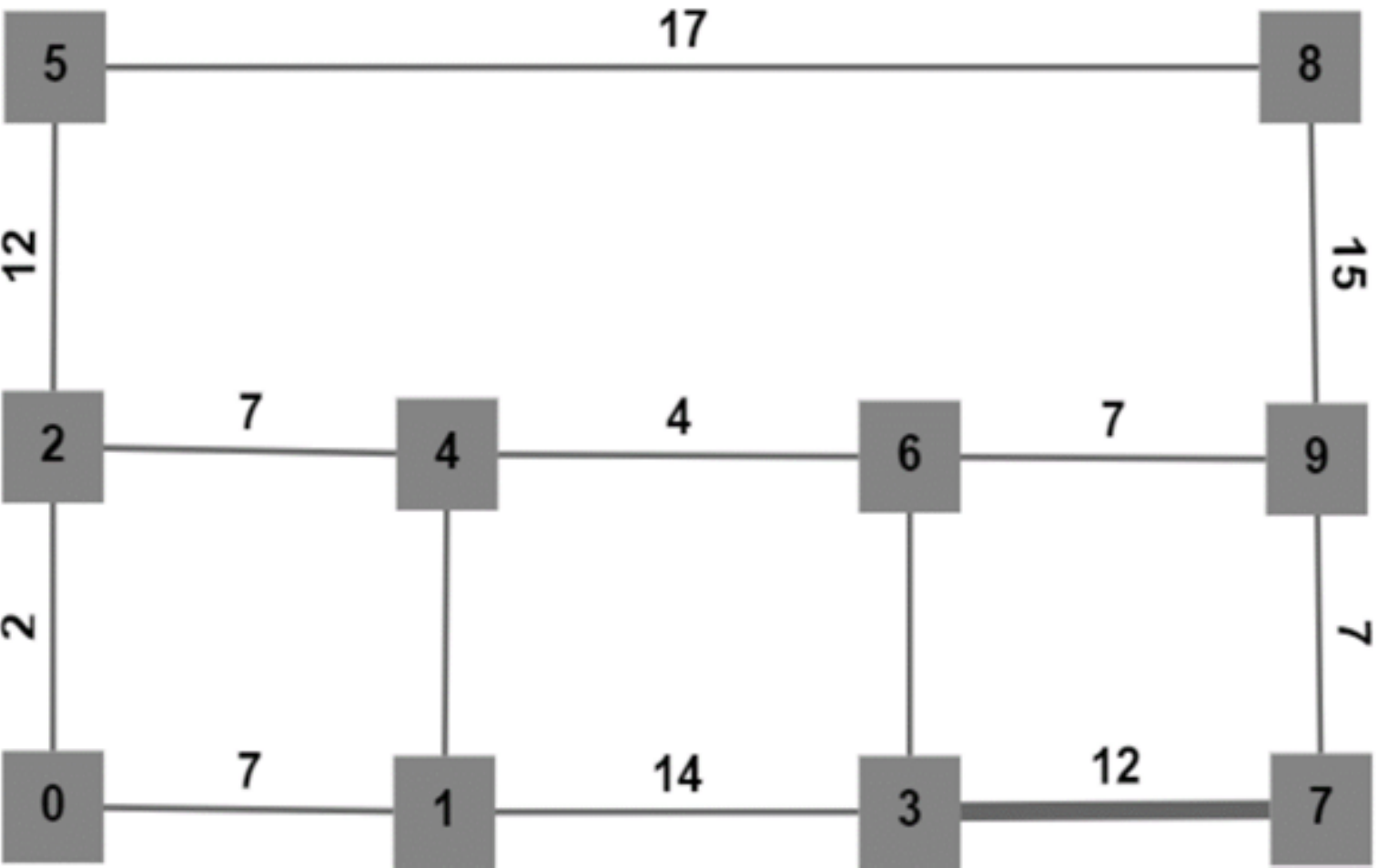} 
			\caption{Road Network at t = 20} 
		\end{minipage} 
		\hspace{0.4cm} 
		\begin{minipage}[t]{4cm} 
			\centering 
			\includegraphics[scale=0.26]{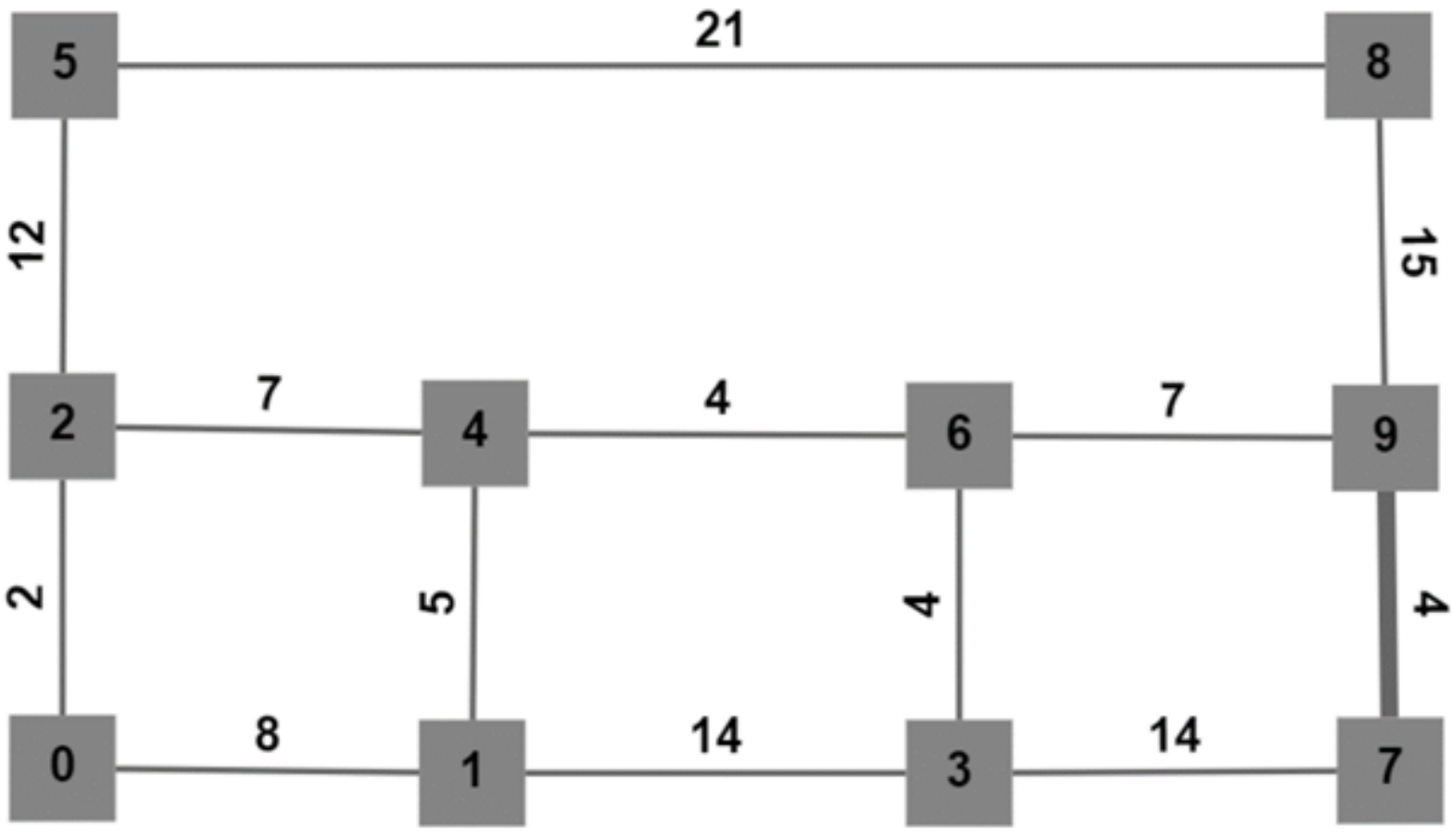} 
			\caption{Road Network at t = 32} 
		\end{minipage} 
	\end{figure}
	
	Table 4 shows the weights of the edges as rows and time stamps as columns. The value in each cell denotes the weight of that edge at a particular time instant. The change from one time stamp to another can be considered as a transition. The eights of those edges that are permissible and are under evaluation at that time stamp are highlighted in bold and underlined. The weights of those edges that constitute the optimal path selected by the proposed algorithm are highlighted in yellow.
	
	%\begin{figure}[htbp]
	%	\centerline{\includegraphics[scale=0.3]{Table3.png}}
	%	\caption{Dynamic Weights of Edges on Transitions}
	%	\label{fig}
	%\end{figure}
	
	% \usepackage[normalem]{ulem}
	\begin{table}[]
		\caption{Dynamic Weights of Edges on Transitions}
		\begin{tabular}{|p{1cm}|c|c|c|c|c|c|c|c|}
			\hline
			 \textbf{Time Stamp/ Edges} & \textbf{t=0} & \textbf{t=6} & \textbf{t=8} & \textbf{t=14} & \textbf{t=18} & \textbf{t=18} & \textbf{t=20} & \textbf{t=32} \\ \hline
			\textbf{E\textsubscript{0-1}}                                                                                        & {\textbf{\hl{ 8 }}}                  & 1                                 & 1                                 & 8                                  & 8                                  & 8                                  & 1                                  & 8                                  \\ \hline
			\textbf{E\textsubscript{0-2}}                                                                                        & {\textbf{\underline{6}}}                  & 3                                 & 2                                 & 3                                  & 2                                  & 2                                  & 2                                  & 2                                  \\ \hline
			\textbf{E\textsubscript{1-3}}                                                                                        & 12                                & 13                                & {\textbf{\hl{ 12 }}}                 & 16                                 & 15                                 & 15                                 & 14                                 & 14                                 \\ \hline
			\textbf{E\textsubscript{1-4}}                                                                                        & 6                                 & 5                                 & {\textbf{\underline{6}}}                  & 5                                  & 5                                  & 5                                  & 5                                  & 5                                  \\ \hline
			\textbf{E\textsubscript{2-4}}                                                                                        & 8                                 & { \textbf{\underline{16}}}                 & 11                                & 11                                 & 1                                  & 1                                  & 7                                  & 7                                  \\ \hline
			\textbf{E\textsubscript{2-5}}                                                                                        & 16                                & { \textbf{\underline{12}}}                 & 12                                & 12                                 & 12                                 & 12                                 & 12                                 & 12                                 \\ \hline
			\textbf{E\textsubscript{3-6}}                                                                                        & 13                                & 4                                 & 4                                 & 4                                  & 4                                  & 4                                  & 4                                  & 4                                  \\ \hline
			\textbf{E\textsubscript{3-7}}                                                                                        & 12                                & 12                                & 12                                & 14                                 & 16                                 & 16                                 & {\textbf{\hl{ 12 }}}                  & 14                                 \\ \hline
			\textbf{E\textsubscript{4-6}}                                                                                        & 6                                 & 6                                 & 4                                 & { \textbf{\underline{4}}}                   & 4                                  & 4                                  & 4                                  & 4                                  \\ \hline
			\textbf{E\textsubscript{5-8}}                                                                                        & 17                                & 17                                & 17                                & 17                                 & { \textbf{\underline{17}}}                  & 17                                 & 17                                 & 21                                 \\ \hline
			\textbf{E\textsubscript{6-9}}                                                                                        & 23                                & 23                                & 23                                & 23                                 & 23                                 & { \textbf{\underline{23}}}                  & 7                                  & 7                                  \\ \hline
			\textbf{E\textsubscript{7-9}}                                                                                        & 16                                & 7                                 & 7                                 & 11                                 & 7                                  & 7                                  & 7                                  & {\textbf{\hl{ 4 }}}                   \\ \hline
			\textbf{E\textsubscript{8-9}}                                                                                        & 17                                & 15                                & 15                                & 15                                 & 15                                 & 15                                 & 15                                 & 15                                 \\ \hline
		\end{tabular}
	\end{table}

	\subsection{Comparison}
	Comparing results from Dijkstra’s algorithm and the routing algorithm proposed in this paper, the final paths selected by both  can be seen in Figures 13 and 14, respectively.
	
	\begin{figure}[htbp]
		\centerline{\includegraphics[scale=0.45]{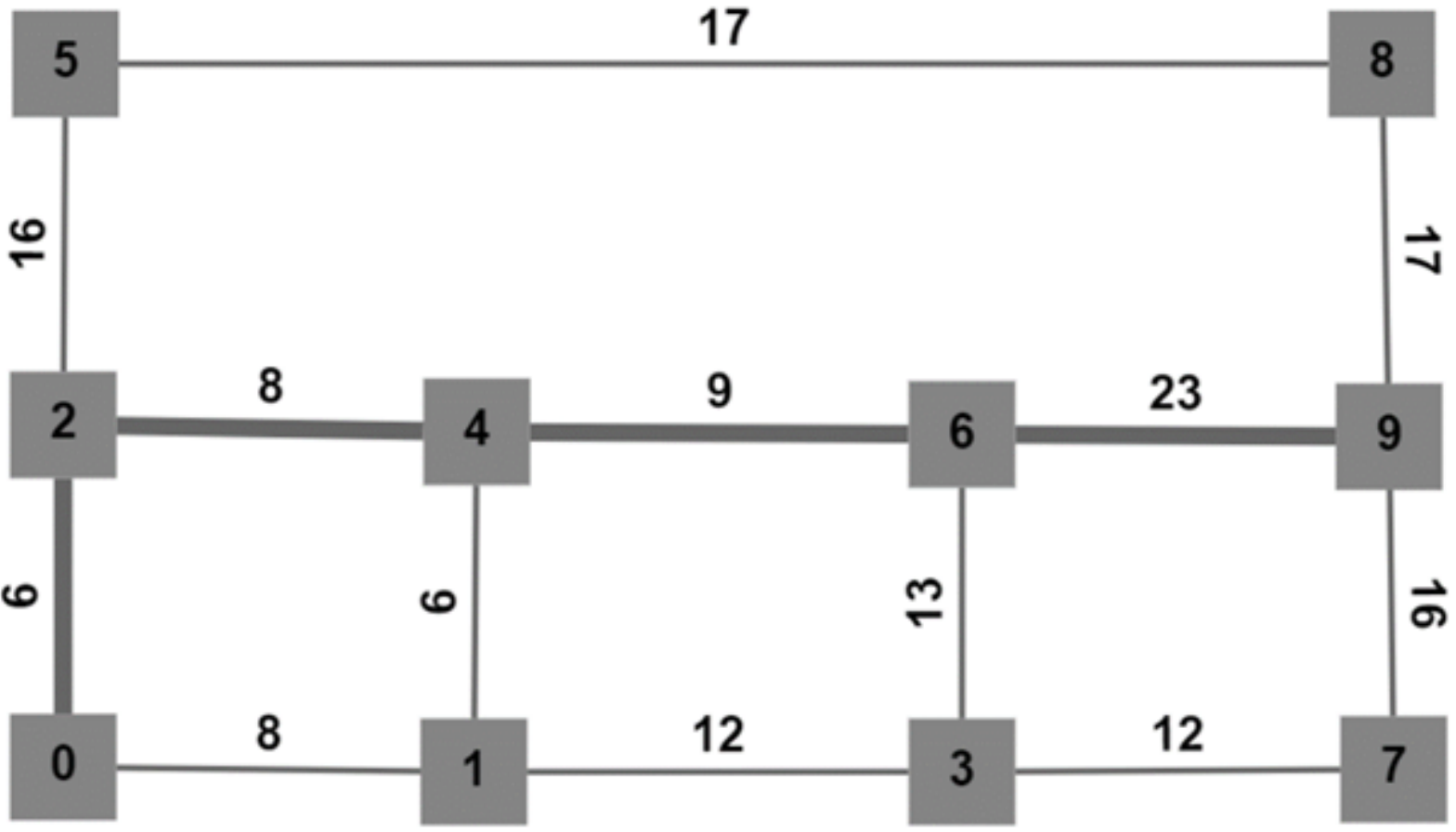}}
		\caption{Shortest path using conventional Dijkstra}
		\label{fig}
	\end{figure}

	\begin{figure}[htbp]
		\centerline{\includegraphics[scale=0.45]{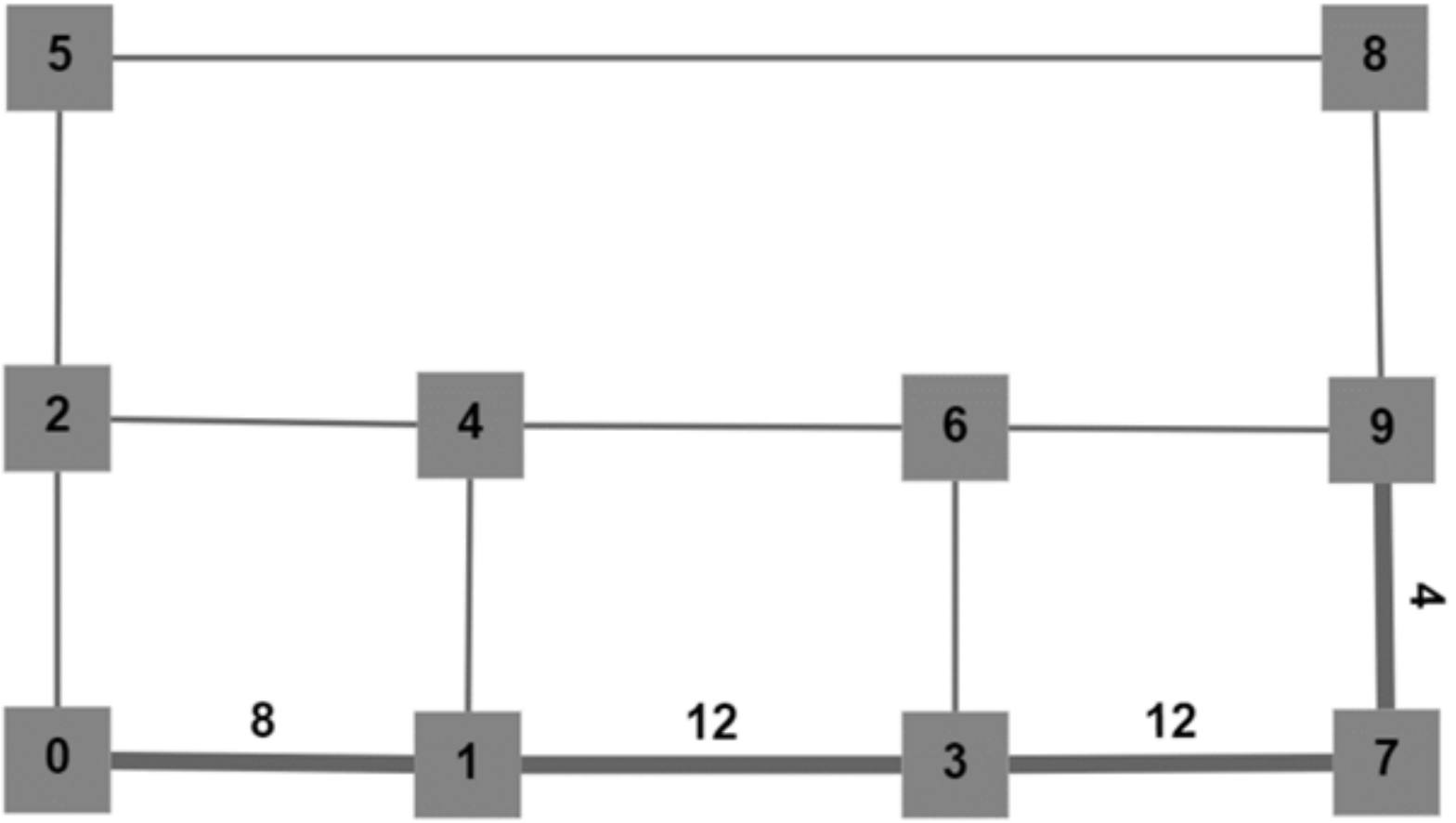}}
		\caption{Shortest path using modified dynamic Dijkstra}
		\label{fig}
	\end{figure}

	A comparative analysis of the above results is carried out next. When using Dijkstra’s algorithm, the shortest path has a total travel time of 46 units. Whereas when using the proposed dynamic routing algorithm with from the traffic prediction model, the shortest path is [ 0 $\rightarrow$ 1 $\rightarrow$ 3 $\rightarrow$ 7 $\rightarrow$ 9 ], with a total travel time of 36 units. This comparison shows that the final path determined by the proposed algorithm is shorter in terms of travel time as compared to the one obtained without considering predicted traffic flows. In this case, the new model leads to a 21.7\% decrease in travel time by simply looking forward in time for traffic movement.

	\subsection{Further Experiments and Results}
	
	As mentioned earlier, the road network for the city of California is used in the form of a graph as an input to the model. From the total of 21,047 nodes and 21,692 edges present in the graph a subset of 80 nodes and 150 edges are chosen as our experimentation set.
	
	In the input data and simulation component of the prediction pipeline, random traffic data is generated over the road network. The model is tested with the specified graph and around 10,000 different $\langle \text{source\_node}, \text{destination\_node} \rangle$ pairs to get the average difference between the least travel time returned by the proposed routing algorithm and conventional Dijkstra’s algorithm. The average difference parameter is given by
	
	\newenvironment{conditions}
	{\par\vspace{\abovedisplayskip}\noindent\begin{tabular}{>{$}l<{$} @{${}={}$} l}}
		{\end{tabular}\par\vspace{\belowdisplayskip}}
	\begin{equation}
		\lambda =  \sum_{i=1}^{N}{\frac{T_i - \tau_i}{N}}
	\end{equation}
	where,
	\begin{conditions}
		$$\lambda$$  &  Average difference \\
		$$\tau_i$$ & Travel time by proposed algorithm for $i^{\textrm{th}}$ test case  \\
		$$T_i$$   & Travel time by Dijkstra's algorithm for $i^{\textrm{th}}$ test case \\
		N &    Total number of test cases   
	\end{conditions}
	
	The average travel time difference is recorded to be 387 seconds, i.e., 6.45 minutes. Model performance was also analysed for different values of $\alpha$ and $\beta$, and the results achieved are shown in Figure 15. This depicts the significant deduction in travel time achieved when using the proposed routing algorithm with traffic prediction as compared to conventional Dijkstra’s algorithm.
	
	Note: In Figure 13, the arbitrariness in results can be attributed to a difference in choice of next edge for the vehicle between simulations for different $\alpha$-$\beta$ values.

	\begin{figure}[htbp]
		\centerline{\includegraphics[scale=0.6]{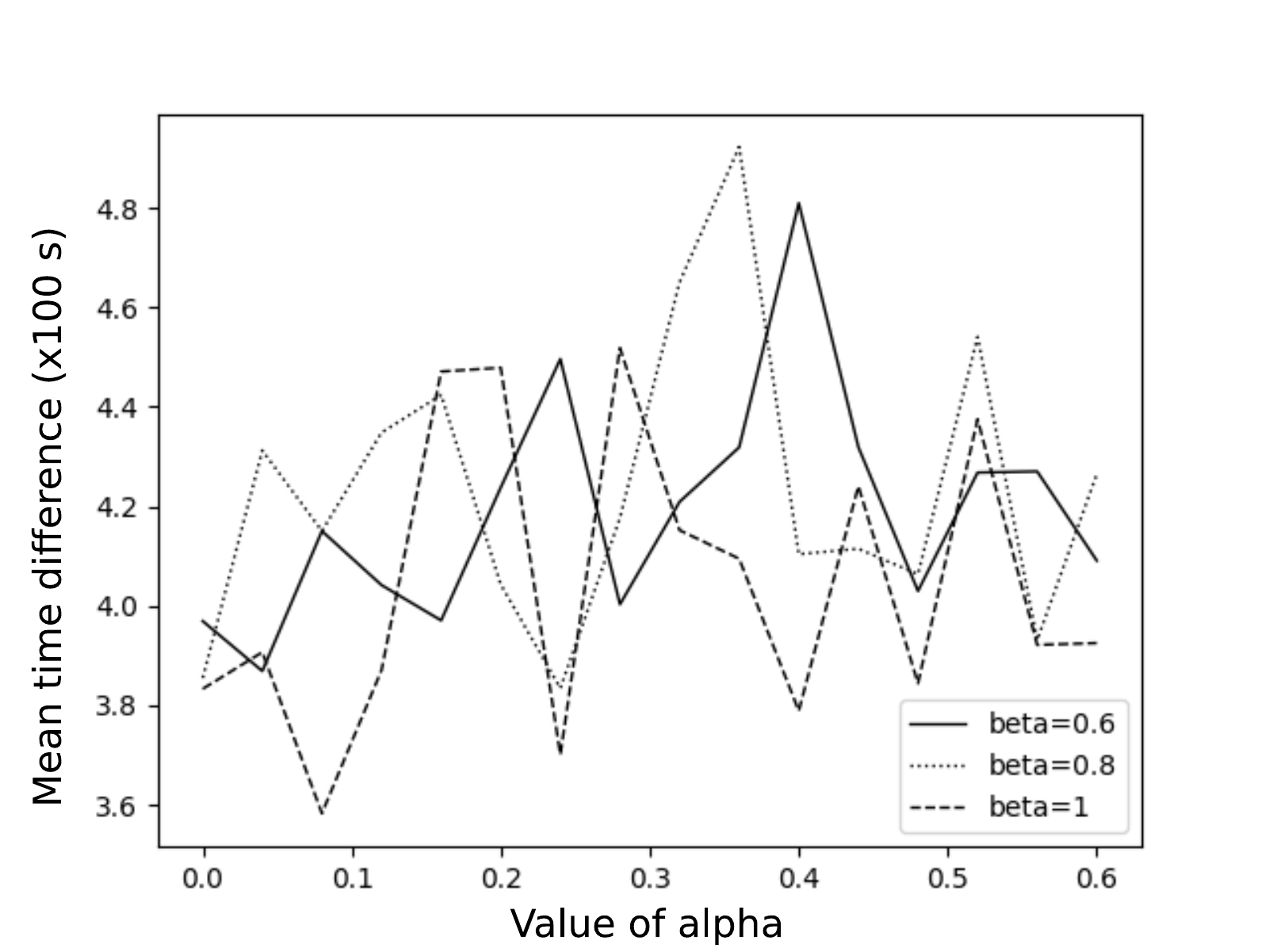}}
		\caption{Mean time difference between the proposed algoritm and Dijkstra's algorithm with varying values of $\alpha$ and $\beta$}
		\label{fig}
	\end{figure}
	
	\section{Conclusion}
	This paper extended Dijkstra’s algorithm to solve the dynamic vehicle routing problem by incorporating traffic movement prediction in the planning stage. It was shown that the weights derived from the traffic flow data dynamically affect the calculations in Dijkstra’s algorithm based on time. The simulation results showed that compared to standard Dijkstra’s algorithm that only considers static traffic data, the proposed routing algorithm with traffic prediction gives optimal routes with a significant reduction in travel time with flowing traffic. This algorithm can be modified in the future for several applications and specializations of the vehicle routing problem by adding constraints pertaining to specific use cases.

\section{Acknowledgement}
The authors would like to thank Mr. Piyush Laddha (Software Product Manager, Jaguar Land Rover India) and Mr. Terence Soares (Software Architect, Jaguar Land Rover India) for their continued support, guidance and encouragement throughout the project. The authors also thank Ms. Isabella Panela (Senior Manager, Jaguar Land Rover Ireland), Mr. Sravan Kallam (Senior Engineer ADAS, Jaguar Land Rover Ireland) and Dr. Damien Dooley (former Senior Engineer ADAS, Jaguar Land Rover Ireland) for reviewing the work and providing valuable feedback for the team to work on. This work was supported by Jaguar Land Rover.

	%\begin{thebibliography}{}
		
		\bibliography{References}
		\bibliographystyle{ieeetr}

\end{document}